\documentclass{amsproc}
\usepackage{amscd,amssymb}

\newtheorem{theorem}{Theorem}[section]
\newtheorem{proposition}[theorem]{Proposition}
\newtheorem{lemma}[theorem]{Lemma}

\theoremstyle{definition}
\newtheorem{definition}[theorem]{Definition}

\theoremstyle{remark}

\numberwithin{equation}{section}

\newcommand{\benu}{\begin{enumerate}\renewcommand{\labelenumi}{{\rm (\roman{enumi})}}\renewcommand{\itemsep}{0pt}}
\newcommand{\eenu}{\end{enumerate}}
\newcommand{\N}{\mathbb{N}}
\newcommand{\Z}{\mathbb{Z}}

\newcommand{\C}{\mathbb{C}}
\newcommand{\T}{\mathbb{T}}

\newcommand{\cK}{{\mathcal K}}
\newcommand{\cL}{{\mathcal L}}

\newcommand{\cO}{{\mathcal O}}

\newcommand{\cM}{{\mathcal M}}

\newcommand{\ip}[2]{\langle{#1},{#2}\rangle}

\DeclareMathOperator{\cspa}{\overline{span}}

\begin{document}
\title{A construction of \boldmath{$C^*$}-algebras 
from \boldmath{$C^*$}-correspondences}
\author{Takeshi Katsura}
\address{
Department of Mathematical Sciences,
University of Tokyo, Komaba, Tokyo 153-8914, JAPAN}
\curraddr{
Department of Mathematics,
University of Oregon, 
Eugene, Oregon 97403-1222, U.S.A.}
\email{katsu@ms.u-tokyo.ac.jp}
\thanks{The author was supported in part by a Research Fellowship 
for Young Scientists of the Japan Society for the Promotion of Science.}

\date{October 31, 2002.}


\begin{abstract}
We introduce a method to define $C^*$-algebras 
from $C^*$-cor\-re\-spon\-denc\-es. 
Our construction generalizes Cuntz-Pimsner algebras, 
crossed products by Hilbert $C^*$-modules, and 
graph algebras. 
\end{abstract}

\maketitle

\setcounter{section}{-1}

\section{Introduction}

In this short article, 
we introduce a new way to define $C^*$-algebras 
from $C^*$-cor\-re\-spon\-denc\-es, 
which generalizes many known constructions in various situations. 
By a $C^*$-cor\-re\-spon\-denc\-e, 
we mean a (right) Hilbert $C^*$-module 
together with a left action (Definition \ref{DefCor}). 
In some articles, this is called a Hilbert $C^*$-bimodule. 
In this paper, the term Hilbert $C^*$-bimodule is reserved to 
denote a right Hilbert $C^*$-module 
which is simultaneously a left Hilbert $C^*$-module 
(Definition \ref{DefBimod}). 
Namely, a Hilbert $C^*$-bimodule has a left inner product, 
but a $C^*$-cor\-re\-spon\-denc\-e may not. 
See Section \ref{ExBimod} in this paper 
for the relation of $C^*$-cor\-re\-spon\-denc\-es 
and Hilbert $C^*$-bimodules. 

In \cite{Pi}, 
Pimsner introduced a class of $C^*$-algebras 
now called Cuntz-Pimsner algebras. 
He constructed $C^*$-algebras from $C^*$-cor\-re\-spon\-denc\-es. 
In \cite{Pi}, 
left actions of $C^*$-cor\-re\-spon\-denc\-es 
are assumed to be injective, 
because it is very possible that 
Cuntz-Pimsner algebras become zero 
if left actions are not injective 
\cite[Remark 1.2 (1)]{Pi}. 
The class of Cuntz-Pimsner algebras is so wide 
that it includes crossed products by automorphisms, 
Cuntz algebras and Cuntz-Krieger algebras, for example. 
However it does not contain 
some interesting examples, such as 
graph algebras of graphs with sinks \cite{FLR} and 
crossed products by partial automorphisms \cite{Ex1}. 
(It was claimed in \cite[Example (4)]{Pi} that 
crossed products by partial automorphisms 
can be obtained by the construction in \cite{Pi}. 
However, it seems that one needs the assumption 
that the ideal $I$ is essential.) 
In \cite{AEE}, 
Abadie, Eilers and Exel gave 
a method to define $C^*$-algebras 
from Hilbert $C^*$-bimodules, 
which generalizes 
crossed products by partial automorphisms. 
Since Hilbert $C^*$-bimodules 
are particular examples of $C^*$-cor\-re\-spon\-denc\-es, 
we can apply the construction in \cite{Pi} 
for Hilbert $C^*$-bimodules. 
This gives the same $C^*$-algebras as in \cite{AEE} 
only when the left actions of the Hilbert $C^*$-bimodules 
are injective. 
For Hilbert $C^*$-bimodules with non-injective left actions, 
the method in \cite{AEE} is different 
from the one in \cite{Pi}. 
In the introduction of \cite{AEE}, 
they mentioned that 
it seems reasonable to expect that a common generalization of 
the constructions in \cite{AEE} and in \cite{Pi} 
could be found. 
Our construction gives this common generalization. 

We construct $C^*$-algebras from arbitrary $C^*$-cor\-re\-spon\-denc\-es. 
We assume no conditions on the $C^*$-cor\-re\-spon\-denc\-es. 
In particular, Hilbert $C^*$-modules need not be full, 
and left actions need be neither injective nor non-degenerate. 
When the left action is injective, 
our construction is same as 
the one of Cuntz-Pimsner algebras. 
When the $C^*$-cor\-re\-spon\-denc\-e comes from 
a Hilbert $C^*$-bimodule, 
our construction is the same as in \cite{AEE}. 
Graph algebras of arbitrary graphs 
can be naturally obtained by our procedure. 
More generally, 
$C^*$-algebras of topological graphs \cite{D,K1} 
are in our class of $C^*$-algebras. 
Thus one of the advantages of our construction 
is that it generalizes many known methods. 
Another advantage is that 
our $C^*$-algebras are closely related to 
the $C^*$-cor\-re\-spon\-denc\-es. 
Namely, the map 
from a $C^*$-cor\-re\-spon\-denc\-e 
to its $C^*$-algebra is always injective, 
and so we can recover all the information about a $C^*$-cor\-re\-spon\-denc\-e 
from its $C^*$-algebra. 
This is not the case for Cuntz-Pimsner algebras 
when the left action of a $C^*$-cor\-re\-spon\-denc\-e 
is not injective. 
Finally, we emphasize that 
even if one starts with 
a $C^*$-cor\-re\-spon\-denc\-e with an injective and non-degenerate 
left action, 
one needs to work with $C^*$-cor\-re\-spon\-denc\-es 
with non-injective or degenerate left actions, 
in order to know about the ideals 
and quotients of its $C^*$-algebra (see \cite{FMR,K2}). 

After preliminaries on Hilbert $C^*$-modules 
and $C^*$-cor\-re\-spon\-denc\-es in Section \ref{SecCorr}, 
we give the definition of our $C^*$-algebras 
constructed from $C^*$-cor\-re\-spon\-denc\-es 
in Section \ref{C*algOfCor}. 
Section \ref{SecEx} is devoted to the study of examples. 

The author is grateful 
to Yasuyuki Kawahigashi for his constant encouragement, 
and to Geoff Price for inviting the author 
to the Conference on Advances in Quantum Dynamics 
and for giving him the opportunity 
to contribute to this proceedings. 
He would also like to thank N. Chris Phillips 
for pointing out many grammatical mistakes 
in the earlier draft. 
This paper was written 
while the author was staying at the University of Oregon. 
He would like to thank people there for their warm hospitality.

\section{\boldmath{$C^*$}-correspondences and their representations.}\label{SecCorr}

\begin{definition}
Let $A$ be a $C^*$-algebra. 
A (right) {\em Hilbert $A$-module} $X$ is a Banach space with
a right action of the $C^*$-algebra $A$ 
and an $A$-valued inner product $\ip{\cdot}{\cdot}_X$ 
satisfying 
\benu
\item $\ip{\xi}{\eta a}_X=\ip{\xi}{\eta}_Xa$, 
\item $\ip{\xi}{\eta}_X=(\ip{\eta}{\xi}_X)^*$, 
\item $\ip{\xi}{\xi}_X\geq 0$ and $\|\xi\|=\|\ip{\xi}{\xi}_X\|^{1/2}$, 
\eenu
for $\xi,\eta\in X$ and $a\in A$. 
\end{definition}

\begin{definition}
For a Hilbert $A$-module $X$, 
we denote by $\cL(X)$ the $C^*$-algebra of all adjointable operators on $X$. 
For $\xi,\eta\in X$, 
the operator $\theta_{\xi,\eta}\in\cL(X)$ is defined 
by $\theta_{\xi,\eta}(\zeta)=\xi\ip{\eta}{\zeta}_X$ for $\zeta\in X$. 
We define $\cK(X)\subset \cL(X)$ by 
$$\cK(X)=\cspa\{\theta_{\xi,\eta}\mid \xi,\eta\in X\},$$
where $\cspa\{\cdots\}$ means the closure of linear span of $\{\cdots\}$. 
\end{definition}

The set $\cK(X)$ is an ideal of $\cL(X)$, 
where an ideal of a $C^*$-algebra always means a two-sided closed ideal, 
which is automatically $*$-invariant. 

\begin{definition}\label{DefCor}
For a $C^*$-algebra $A$, 
we say that $X$ is a {\em $C^*$-cor\-re\-spon\-denc\-e} over $A$ 
when $X$ is a Hilbert $A$-module and 
a $*$-homomorphism $\varphi_X:A\to \cL(X)$ is given. 
\end{definition}

We refer to $\varphi_X$ as the left action of 
a $C^*$-cor\-re\-spon\-denc\-e $X$. 

\begin{definition}
A $C^*$-cor\-re\-spon\-denc\-e $X$ over a $C^*$-algebra $A$ is called 
{\em faithful} if $\varphi_X$ is injective, 
{\em non-degenerate} 
if $\cspa\{\varphi_X(a)\xi\in X\mid a\in A,\xi\in X\}=X$, 
and {\em full} if 
$\cspa\{\ip{\xi}{\eta}_X\in A\mid \xi,\eta\in X\}=A$. 
\end{definition}

\begin{definition}
A {\em representation} of a $C^*$-cor\-re\-spon\-denc\-e $X$ over $A$ 
on a $C^*$-algebra $B$ is 
a pair $(\pi,t)$ consisting of 
a $*$-homomorphism $\pi:A\to B$ and 
a linear map $t:X\to B$ satisfying 
\benu
\item $t(\xi)^*t(\eta)=\pi\big(\ip{\xi}{\eta}_X\big)$ 
for $\xi,\eta\in X$,
\item $\pi(a)t(\xi)=
t\big(\varphi_X(a)\xi\big)$ for $a\in A$, $\xi\in X$. 
\eenu
\end{definition}

A representation of a $C^*$-cor\-re\-spon\-denc\-e is called 
an isometric covariant representation in \cite{MS}. 
For a representation $(\pi,t)$ on $B$, 
we denote by $C^*(\pi,t)$ the $C^*$-algebra generated 
by the images of $\pi$ and $t$ in $B$. 
Note that for a representation $(\pi,t)$ of $X$, 
we have $t(\xi)\pi(a)=t(\xi a)$ automatically 
because the condition (i) above, 
combined with the fact that $\pi$ is a $*$-homomorphism, implies 
$$\big\|t(\xi)\pi(a)-t(\xi a)\big\|^2
=\big\|\big(t(\xi)\pi(a)-t(\xi a)\big)^*
 \big(t(\xi)\pi(a)-t(\xi a)\big)\big\|
=0.$$
Note also that for $\xi\in X$, 
we have $\|t(\xi)\|\leq\|\xi\|$ 
because 
$$\|t(\xi)\|^2=\|t(\xi)^*t(\xi)\|=\|\pi(\ip{\xi}{\xi}_X)\|
\leq \|\ip{\xi}{\xi}_X\|=\|\xi\|^2.$$

\begin{definition}
A representation $(\pi,t)$ is said to be {\em injective}
when the $*$-homomorphism $\pi$ is injective. 
\end{definition}

By the above computation, 
we see that $t$ is isometric 
for an injective representation $(\pi,t)$. 

\begin{definition}
For a representation $(\pi,t)$ of a $C^*$-cor\-re\-spon\-denc\-e $X$ on $B$, 
we define a $*$-homomorphism $\psi_t:\cK(X)\to B$ 
by $\psi_t(\theta_{\xi,\eta})=t(\xi)t(\eta)^*\in B$ 
for $\xi,\eta\in X$. 
\end{definition}

For the well-definedness of the $*$-homomorphism $\psi_t$, 
see \cite[Lemma 2.2]{KPW}, for example.
Note that $\psi_t\big(\cK(X)\big)\subset C^*(\pi,t)$. 
Note also that 
$\psi_t$ is injective for an injective representation $(\pi,t)$.

\section{\boldmath{$C^*$}-algebras associated with \boldmath{$C^*$}-cor\-re\-spon\-denc\-es}\label{C*algOfCor}

\begin{definition}
For an ideal $I$ of a $C^*$-algebra $A$, 
we define $I^{\perp}\subset A$ by 
$$I^{\perp}=\{a\in A\mid ab=0 \mbox{ for all }b\in I\}.$$
\end{definition}

Recall that an ideal $I$ is called {\em essential} if $I^{\perp}=0$. 

\begin{lemma}\label{Iperp}
For an ideal $I$ of a $C^*$-algebra $A$, 
$I^{\perp}$ is an ideal of $A$, 
and for an ideal $J$ of $A$, 
we have $J\subset I^{\perp}$ if and only if $J\cap I=\{0\}$. 
\end{lemma}

\begin{proof}
It is clear that $I^{\perp}$ is an ideal. 
It is also clear that $I^{\perp}\cap I=\{0\}$. 
Hence if an ideal $J$ satisfies $J\subset I^{\perp}$, 
then we have $J\cap I=\{0\}$. 
Conversely, suppose that an ideal $J$ satisfies $J\cap I=\{0\}$. 
Since $ab\in J\cap I=\{0\}$ for any $a\in J$ and $b\in I$, 
we have $J\subset I^{\perp}$. 
This completes the proof. 
\end{proof}

Lemma \ref{Iperp} means that 
$I^{\perp}$ is the maximal ideal satisfying $I^{\perp}\cap I=\{0\}$. 

\begin{definition}
For a $C^*$-cor\-re\-spon\-denc\-e $X$ over $A$ , 
we define an ideal $J_X$ of $A$ by 
$$J_X=\varphi_X^{-1}\big(\cK(X)\big)\cap \big(\ker\varphi_X\big)^{\perp}.$$
\end{definition}

Note that $J_X=\varphi_X^{-1}\big(\cK(X)\big)$ 
when $\varphi_X$ is injective. 
By Lemma \ref{Iperp}, 
we see that the ideal $J_X$ is the maximal ideal 
on which the restriction of $\varphi_X$ is an injection into $\cK(X)$. 
From this fact, we easily get the following. 

\begin{lemma}\label{J=JX}
If there exists an ideal $J$ of $A$ such that 
the restriction of $\varphi_X$ to $J$ 
is an isomorphism onto $\cK(X)$, 
then $J=J_X$. 
\end{lemma}

\begin{definition}
A representation $(\pi,t)$ is said to satisfy 
{\em Condition $(*)$} 
if we have $\pi(a)=\psi_t(\varphi_X(a))$ 
for all $a\in J_X$. 
\end{definition}

\begin{definition}\label{DefOX}
For a $C^*$-cor\-re\-spon\-denc\-e $X$ over a $C^*$-algebra $A$, 
the $C^*$-algebra $\cO_X$ is defined by $\cO_X=C^*(\pi_X,t_X)$ 
where $(\pi_X,t_X)$ is the universal representation of $X$ 
satisfying Condition $(*)$. 
\end{definition}

One can easily show the existence of the universal representation of $X$ 
satisfying Condition $(*)$ (see, for example, \cite{Bl}). 
One can also define $\cO_X$ in a more concrete way 
using a Fock space as in \cite{Pi} or \cite{MS}. 
We have that $\cO_X$ is isomorphic to 
the relative Cuntz-Pimsner $C^*$-algebra $\cO(J_X,X)$ 
defined in \cite[Definition 2.18]{MS} 
(see Example \ref{relCP} in this paper). 
Since $J_X\cap\ker\varphi_X=\{0\}$, 
\cite[Proposition 2.21]{MS} gives us the following. 

\begin{proposition}\label{isometric}
The universal representation $(\pi_X,t_X)$ of $X$ 
on $\cO_X$ is injective.
\end{proposition}

In \cite{K2}, 
we show that the $C^*$-algebra $\cO_X$ 
is the smallest $C^*$-algebra 
among $C^*$-algebras $C^*(\pi,t)$ 
generated by representations $(\pi,t)$ of $X$ 
which are injective and admit certain actions of $\T$ 
called gauge actions. 
This fact gives another definition of $\cO_X$ 
without using Condition $(*)$ or $J_X$.

\section{Examples}\label{SecEx}

\subsection{Cuntz-Pimsner algebras}

In \cite{Pi}, Pimsner defined two kinds of $C^*$-algebras 
from a $C^*$-cor\-re\-spon\-denc\-e $X$, 
which are now called the Toeplitz algebra and 
the Cuntz-Pimsner algebra of $X$, respectively.
He also defined the augmented ones. 
When the $C^*$-cor\-re\-spon\-denc\-e $X$ is faithful, 
the $C^*$-algebra $\cO_X$ defined in this paper
is the same as 
the augmented Cuntz-Pimsner algebra of $X$. 
When the $C^*$-cor\-re\-spon\-denc\-e $X$ is full, 
the augmented Cuntz-Pimsner algebra of $X$ 
coincides with its 
Cuntz-Pimsner algebra.
Hence our $C^*$-algebra $\cO_X$ is the same as 
the Cuntz-Pimsner algebra of $X$ 
when the $C^*$-cor\-re\-spon\-denc\-e $X$ is both faithful and full. 
For a general $C^*$-cor\-re\-spon\-denc\-e $X$, 
we have a surjection from our $C^*$-algebra $\cO_X$ 
onto the augmented Cuntz-Pimsner algebra of $X$. 
This surjection is never injective 
unless the $C^*$-cor\-re\-spon\-denc\-e $X$ is faithful.

\subsection{Relative Cuntz-Pimsner algebras}\label{relCP}

Let $X$ be a $C^*$-cor\-re\-spon\-denc\-e over $A$. 
Let $J_0$ be an ideal of $A$ such that $\varphi_X(J_0)\subset\cK(X)$. 
The {\em relative Cuntz-Pimsner algebra} $\cO(J_0,X)$ is 
the universal $C^*$-algebra generated by representations $(\pi,t)$ 
of $X$ satisfying $\pi(a)=\psi_t(\varphi_X(a))$ for all $a\in J_0$. 
Actually Muhly and Solel defined relative Cuntz-Pimsner algebras 
in more concrete way using Fock spaces in \cite[Definition 2.18]{MS}, 
and then proved that they have the universal property 
in \cite[Theorem 2.19]{MS}. 
Note that when $J_0=0$, $\cO(J_0,X)$ 
is the augmented Toeplitz algebra of $X$, 
and when $J_0=\varphi_X^{-1}\big(\cK(X)\big)$, 
$\cO(J_0,X)$ 
is the augmented Cuntz-Pimsner algebra. 
Our $C^*$-algebra $\cO_X$ is isomorphic to $\cO(J_X,X)$. 
In \cite{K2}, 
we prove that 
the relative Cuntz-Pimsner algebra $\cO(J_0,X)$ 
is isomorphic to $\cO_{X'}$ for 
a certain $C^*$-cor\-re\-spon\-denc\-e $X'$. 
This observation helps us to study 
relative Cuntz-Pimsner algebras.

\subsection{Crossed products by Hilbert \boldmath{$C^*$}-bimodules}\label{ExBimod}

\begin{definition}\label{DefBimod}
Let $A$ be a $C^*$-algebra. 
A {\em Hilbert $A$-bimodule} $X$ is 
a Hilbert $A$-module 
together with a left action $\varphi_X:A\to\cL(X)$ and 
a left inner product ${}_X\ip{\cdot}{\cdot}:X\times X\to A$, 
which satisfy 
\begin{enumerate}
\renewcommand{\labelenumi}{{\rm (\roman{enumi})$'$}}
\renewcommand{\itemsep}{0pt}
\item ${}_X\ip{\varphi_X(a)\xi}{\eta}=a\cdot{}_X\ip{\xi}{\eta}$, 
\item ${}_X\ip{\xi}{\eta}={}_X\ip{\eta}{\xi}^*$, 
\item ${}_X\ip{\xi}{\xi}\geq 0$, 
\end{enumerate}
for $\xi,\eta\in X$, $a\in A$ and 
\begin{equation}\label{HilBimod}
\varphi_X\big({}_X\ip{\xi}{\eta}\big)\zeta=\xi\ip{\eta}{\zeta}_X,
\end{equation}
for $\xi,\eta,\zeta\in X$. 
\end{definition}

In the original definition (\cite[Defnition 1.8]{BMS}), 
it was not assumed that $\varphi_X(a):X\to X$ is adjointable 
and that its adjoint is $\varphi_X(a^*)$. 
However this fact follows from (\ref{HilBimod})
(see \cite{BMS}). 
Following \cite{BMS}, 
we define an ideal $I_X\subset A$ by
$$I_X=\cspa\{{}_X\ip{\xi}{\eta}\in A\mid \xi,\eta\in X\}.$$ 

\begin{lemma}\label{aIX=0}
Let $X$ be a Hilbert $A$-bimodule. 
For $a\in A$, we have $\varphi_X(a)=0$ 
if and only if $a\in I_X^{\perp}$. 
\end{lemma}

\begin{proof}
If $\varphi_X(a)=0$ 
then we have $a\cdot{}_X\ip{\xi}{\eta}={}_X\ip{\varphi_X(a)\xi}{\eta}=0$ 
for all $\xi,\eta\in X$. 
Hence $a\in I_X^{\perp}$. 
If $a\in I_X^{\perp}$ 
then we have 
$${}_X\ip{\varphi_X(a)\xi}{\varphi_X(a)\xi}=a\cdot{}_X\ip{\xi}{\varphi_X(a)\xi}=0$$
for all $\xi\in X$. 
Hence $\varphi_X(a)=0$. 
\end{proof}

The equation (\ref{HilBimod}) means that 
$$\varphi_X({}_X\ip{\xi}{\eta})
=\theta_{\xi,\eta},$$
for $\xi,\eta\in X$. 
Hence the restriction of $\varphi_X$ to the ideal $I_X$ 
is a surjection onto $\cK(X)$. 
By Lemma \ref{aIX=0}, this restriction is injective. 
Thus the restriction of $\varphi_X$ to the ideal $I_X$ 
is an isomorphism onto $\cK(X)$ (\cite[Proposition 1.10]{BMS}). 
Hence Lemma \ref{J=JX} gives us the following. 

\begin{lemma}\label{Bimod=>Corr}
Let $A$ be a $C^*$-algebra, 
and $X$ be a Hilbert $A$-bimodule. 
If we consider $X$ as a $C^*$-cor\-re\-spon\-denc\-e over $A$, 
then we have $J_X=I_X$. 
\end{lemma}

By Lemma \ref{Bimod=>Corr}, 
the left inner product ${}_X\ip{\cdot}{\cdot}$ of 
a Hilbert $A$-bimodule $X$ is 
completely determined by the structure of $X$ 
as a $C^*$-cor\-re\-spon\-denc\-e.
Namely, for $\xi,\eta\in X$, 
${}_X\ip{\xi}{\eta}\in A$ is the unique element $a\in J_X$ 
with $\varphi_X(a)=\theta_{\xi,\eta}$. 
We can reverse this procedure. 

\begin{lemma}\label{Corr=>Bimod}
Let $A$ be a $C^*$-algebra, 
and $X$ be a $C^*$-cor\-re\-spon\-denc\-e over $A$. 
If $\varphi_X(J_X)=\cK(X)$, 
then $X$ is a Hilbert $A$-bimodule 
with the left inner product ${}_X\ip{\cdot}{\cdot}$ given by 
$${}_X\ip{\xi}{\eta}
=\big(\varphi_X|_{J_X}\big)^{-1}(\theta_{\xi,\eta})\in J_X\subset A$$
for $\xi,\eta\in X$. 
\end{lemma}

\begin{proof}
Straightforward. 
\end{proof}

By the above argument, 
we can say that Hilbert $A$-bimodules are nothing but 
$C^*$-cor\-re\-spon\-denc\-es $X$ over $A$ satisfying $\varphi_X(J_X)=\cK(X)$, 
and left inner products of Hilbert $A$-bimodules 
are uniquely determined by the structures 
as $C^*$-cor\-re\-spon\-denc\-es. 
A $C^*$-cor\-re\-spon\-denc\-e arising from a Hilbert $C^*$-bimodule 
is always non-degenerate, 
and it is faithful if and only if 
$J_X$ is an essential ideal in $A$ 
by Lemma \ref{aIX=0} and Lemma \ref{Bimod=>Corr}. 
In \cite{AEE}, 
Abadie, Eilers and Exel defined 
the crossed product $A\rtimes_X \Z$ by a Hilbert $A$-bimodule $X$ 
as follows. 

\begin{definition}[{\cite[Definition 2.1]{AEE}}]\label{DefCovRep}
Let $X$ be a Hilbert $A$-bimodule. 
A {\em covariant representation} of $X$ 
on a $C^*$-algebra $B$ is 
a pair $(\pi,t)$ consisting of 
a $*$-homomorphism $\pi:A\to B$ 
and a linear map $t:X\to B$ satisfying 
\benu
\item $t(\varphi_X(a)\xi)=\pi(a)t(\xi)$, 
\item $t(\xi a)=t(\xi)\pi(a)$, 
\item $\pi({}_X\ip{\xi}{\eta})=t(\xi)t(\eta)^*$, 
\item $\pi(\ip{\xi}{\eta}_X)=t(\xi)^*t(\eta)$, 
\eenu
for $a\in A$ and $\xi,\eta\in X$. 
\end{definition}

In \cite{AEE}, 
they considered only the case that $B=B(H)$ 
for some Hilbert space $H$. 

\begin{definition}[{\cite[Definition 2.4]{AEE}}]
Let $A$ be a $C^*$-algebra, 
and $X$ be a Hilbert $A$-bimodule. 
The {\em crossed product} $A\rtimes_X\Z$ of $A$ by $X$ 
is the universal $C^*$-algebra generated 
by the images of covariant representations of $X$. 
\end{definition}

The conditions (i) and (iv) in Definition \ref{DefCovRep} are the same 
as the conditions of representations of $X$ as a $C^*$-cor\-re\-spon\-denc\-e. 
The condition (ii) follows from (iv) 
(and also the condition (i) follows from (iii), 
but we do not use this fact). 
Since $J_X=\cspa\{{}_X\ip{\xi}{\eta}\in A\mid \xi,\eta\in X\}$ 
and $\varphi_X({}_X\ip{\xi}{\eta})=\theta_{\xi,\eta}$ 
for $\xi,\eta\in X$, 
we see that the condition (iii) is 
equivalent to Condition $(*)$. 
Thus we get the following. 

\begin{proposition}
Let $A$ be a $C^*$-algebra, 
and $X$ be a Hilbert $A$-bimodule. 
The crossed product $A\rtimes_X\Z$ of $A$ by $X$ 
is canonically isomorphic to $\cO_X$ 
where $X$ is considered as a $C^*$-cor\-re\-spon\-denc\-e over $A$. 
\end{proposition}

\subsection{Graph algebras}\label{SecGra}

\begin{definition}\label{Defgrph}
A (directed) {\em graph} $E=(E^0,E^1,r,s)$ consists 
of two sets $E^0,E^1$ which are sets of vertices and edges 
respectively, 
and two maps $r,s:E^1\to E^0$ which indicate 
the range $r(e)$ and the source $s(e)$ of a directed edge $e\in E^1$. 
\end{definition}

In \cite{FLR}, 
Fowler, Laca, and Raeburn defined graph algebras 
for arbitrary graphs, 
generalizing the definitions of \cite{KPRR} and \cite{KPR}. 

\begin{definition}
The {\em graph algebra} $C^*(E)$ of a graph $E$ 
is the universal $C^*$-algebra generated 
by mutually orthogonal projections $\{p_v\}_{v\in E^0}$ 
and partial isometries $\{s_e\}_{e\in E^1}$ 
with orthogonal ranges, 
such that $s_e^*s_e=p_{r(e)}$, $s_es_e^*\leq p_{s(e)}$ 
for $e\in E^1$ and 
\begin{equation}\label{CPFam}
p_v=\sum_{e\in s^{-1}(v)}s_es_e^*\quad \mbox{ if }0<|s^{-1}(v)|<\infty.
\end{equation}
\end{definition}

Take a graph $E=(E^0,E^1,r,s)$. 
Set $A=C_0(E^0)$. 
The linear space $C_c(E^1)$ is a pre-Hilbert $A$-module 
under the following operations 
\begin{align*}
\ip{\xi}{\eta}_X(v)=\sum_{e\in r^{-1}(v)}\overline{\xi(e)}\eta(e)\in\C
\quad \mbox{for }v\in E^0,\\
(\xi f)(e)=\xi(e)f(r(e))\in\C
\quad \mbox{for }e\in E^1,
\end{align*}
for $\xi,\eta\in C_c(E^1)$ and $f\in A$. 
The completion $X=X(E)$ of $C_c(E^1)$ 
in the norm defined by $\|\xi\|=\|\ip{\xi}{\xi}_X\|^{1/2}$ 
is a Hilbert $A$-module. 
Define $\varphi_X:A\to\cL(X)$ by 
$$\varphi_X(f)\xi(e)=f(s(e))\xi(e)
\quad \mbox{for }e\in E^1,$$
for $f\in A$ and $\xi\in C_c(E^1)\subset X$. 
We have a $C^*$-cor\-re\-spon\-denc\-e $X$ over $A$. 
We see that 
\begin{align*}
\varphi_X^{-1}\big(\cK(X)\big)
&=C_0\big(\{v\in E^0\mid |s^{-1}(v)|<\infty\}\big),\\
\ker\varphi_X &= C_0\big(\{v\in E^0\mid |s^{-1}(v)|=0\}\big). 
\end{align*}
Hence we get $J_X=C_0\big(\{v\in E^0\mid 0<|s^{-1}(v)|<\infty\}\big)$. 

If we have mutually orthogonal projections $\{p_v\}_{v\in E^0}$ 
and partial isometries $\{s_e\}_{e\in E^1}$ 
with orthogonal ranges, 
then we can define a $*$-homomorphism $\pi$ 
and a linear map $t$ 
by 
\begin{align*}
\pi(f)&=\sum_{v\in E^0}f(v)p_v\quad \mbox{ for }f\in C_0(E^0)=A,\\ 
t(\xi)&=\sum_{e\in E^1}\xi(e)s_e\quad \mbox{ for }
\xi\in C_c(E^1)\subset X(E). 
\end{align*}
This pair $(\pi,t)$ is a representation of $X(E)$ 
if and only if 
two conditions $s_e^*s_e=p_{r(e)}$, $s_es_e^*\leq p_{s(e)}$ 
are satisfied for $e\in E^1$. 
This representation $(\pi,t)$ satisfies Condition $(*)$ 
whenever (\ref{CPFam}) is fulfilled. 
Hence we get the following. 

\begin{proposition}
The graph algebra $C^*(E)$ of a graph $E$ 
is isomorphic to the $C^*$-algebra $\cO_{X(E)}$ 
of the $C^*$-cor\-re\-spon\-denc\-e $X(E)$. 
\end{proposition}

A vertex $v\in E^0$ is called a sink if $s^{-1}(v)=\emptyset$, 
and a source if $r^{-1}(v)=\emptyset$. 
The $C^*$-cor\-re\-spon\-denc\-e $X$ defined by a graph $E$ 
is faithful if and only if $E$ has no sinks, 
and full if and only if $E$ has no sources. 
In particular, 
for a graph $E$ with no sinks, 
the graph algebra $C^*(E)$ 
is isomorphic to 
the augmented Cuntz-Pimsner algebra $\widetilde{\cO}_X$ 
of the $C^*$-cor\-re\-spon\-denc\-e $X$ defined by the graph $E$ 
(\cite[Proposition 12]{FLR}). 

\subsection{\boldmath{$C^*$}-algebras arising from topological graphs}\label{SecTopGra}

In \cite{K1}, 
we generalize the construction of graph algebras 
to topological graphs (see also \cite{D}). 

\begin{definition}[{\cite[Definition 2.1]{K1}}]
A {\em topological graph} 
is a quadruple $E=(E^0,E^1,d,r)$ 
where $E^0,E^1$ are locally compact spaces, 
$d:E^1\to E^0$ is a local homeomorphism 
and $r:E^1\to E^0$ is a continuous map. 
\end{definition}

The maps $d,r:E^1\to E^0$ indicate 
the domain and source map, respectively.
For a topological graph $E=(E^0,E^1,d,r)$, 
the triple $(E^1,d,r)$ is called 
a topological correspondence over $E^0$ 
in \cite{K1}. 
From a topological graph $E=(E^0,E^1,d,r)$, 
we can define a $C^*$-cor\-re\-spon\-denc\-e $C_d(E^1)$ over $C_0(E^0)$ 
by 
$$C_d(E^1)=\{\xi\in C(E^1)\mid \ip{\xi}{\xi}\in C_0(E^0)\},$$
where the inner product, the right action, 
and the left action are defined 
in the same way as in Example \ref{SecGra}, 
but we use the domain map $d$ 
to define a Hilbert $C_0(E^0)$-module, 
and the range map $r$ to define a left action. 
This is opposite from the case of graph algebras. 
The author believes that this convention is more natural 
than the one used in the theory of graph algebras, 
and this convention behaves well 
when we consider topological graphs 
as a kind of dynamical system (see Example \ref{ExParMor}). 
We can show that $C_c(E^1)$ is dense in $C_d(E^1)$. 
We have $J_{C_d(E^1)}=C_0(E^0_{\mathrm{rg}})$ 
where 
\begin{align*}
E^0_{\mathrm{rg}}
=\big\{v\in E^0\ \big|\ &\mbox{ there exists a neighborhood } 
V \mbox{ of } v\\
&\mbox{ such that }
r^{-1}(V)\subset E^1 \mbox{ is compact, and }r(r^{-1}(V))=V\big\}.
\end{align*}
We define the $C^*$-algebra $\cO(E)$ of 
a topological graph $E$ 
in the same way as we define the $C^*$-algebra $\cO_{X}$ 
in this paper for the $C^*$-cor\-re\-spon\-denc\-e $X=C_d(E^1)$ 
\cite[Definition 2.10]{K1}. 
The $C^*$-cor\-re\-spon\-denc\-e $C_d(E^1)$ over $C_0(E^0)$ 
defined from a topological graph $E$ 
is always non-degenerate. 
It is full if and only if $d$ is surjective, 
and faithful if and only if the image of $r$ is dense in $E^0$. 

Note that we can define $C^*$-cor\-re\-spon\-denc\-es 
from so-called continuous measured graphs 
which are a generalization of topological graphs 
(see \cite[Example 1.7]{Sc1}).

\subsection{Crossed products by partial morphisms}\label{ExParMor}

Let $X,Y$ be locally compact spaces, 
and let $\sigma:X\to Y$ be a continuous map. 
The map $\sigma$ determines a $*$-homomorphism 
$\varphi:C_0(Y)\to C_b(X)$ by $\varphi(f)(x)=f(\sigma(x))$. 
The $C^*$-algebra $C_b(X)$ 
is isomorphic to the multiplier algebra of $C_0(X)$, 
and the map $\varphi$ is non-degenerate in the sense that 
$\cspa\{\varphi(f)g\in C_0(X)\mid f\in C_0(Y),g\in C_0(X)\}=C_0(X)$. 
Every non-degenerate $*$-homomorphism 
from $C_0(Y)$ to $C_b(X)$ is obtained by this procedure. 
Thus the non-commutative analogues of continuous maps between 
locally compact spaces are 
non-degenerate $*$-homomorphisms from one $C^*$-algebra $A$ 
to the multiplier algebra $\cM(B)$ of another $C^*$-algebra $B$. 
Such a map is called a {\em morphism} from $A$ to $B$ in \cite{Lnc}. 
On the other hand, 
a $*$-homomorphism $\varphi:C_0(Y)\to C_0(X)$ corresponds 
to a proper continuous map from some open subset of $X$ to $Y$. 
This open subset is determined by the hereditary subalgebra 
in $C_0(X)$ generated by the image of $\varphi$. 
Now it is natural to consider 
a continuous map from some open subset of $X$ to $Y$ 
and its non-commutative analogue.

\begin{definition}[{Cf. \cite[Example 1.2]{Sc1}}]
Let $A,B$ be $C^*$-algebras. 
A {\em partial morphism} from $A$ to $B$ 
is a non-degenerate $*$-homomorphism $\varphi$ from $A$ to 
the multiplier algebra $\cM(B_0)$ 
of some hereditary subalgebra $B_0$ of $B$. 
\end{definition}

Partial morphisms between commutative $C^*$-algebras 
correspond bijectively to partially defined continuous maps 
between their spectra. 
Any $*$-homomorphism from $A$ to $B$ determines 
a partial morphism from $A$ to $B$. 
More generally, we have the following.

\begin{definition}
Let $A$, $B$ be $C^*$-algebras, 
let $I$ be an ideal of $A$, and 
let $B'$ be a $C^*$-subalgebra of $B$. 
Let $\varphi':I\to \cM(B')$ be a $*$-homomorphism. 
Define a hereditary subalgebra $B_0$ of $B$ by 
$$B_0=\cspa\big\{\varphi'(x_1)b_1'bb_2'\varphi'(x_2)\in B\ \big|\
x_1,x_2\in I, b_1',b_2'\in B', b\in B\big\}.$$
We define a $*$-homomorphism $\varphi:A\to \cM(B_0)$ 
by 
$$\varphi(a)b_0=\varphi'(ax_1)b_1'bb_2'\varphi'(x_2),\quad 
b_0\varphi(a)=\varphi'(x_1)b_1'bb_2'\varphi'(x_2a)\in B_0,$$ 
for $a\in A$ and 
$b_0=\varphi'(x_1)b_1'bb_2'\varphi'(x_2)\in B_0$. 
It is standard to see that 
this gives a well-defined non-degenerate 
$*$-homomorphism $\varphi:A\to\cM(B_0)$. 
We call this partial morphism $\varphi:A\to\cM(B_0)$ 
the partial morphism determined by $\varphi':I\to \cM(B')$. 
\end{definition}

\begin{definition}[{\cite[Definition 3.1]{Ex1}}]
A {\em partial automorphism} of a $C^*$-algebra $A$ 
is a triple $(\theta,I,J)$ where $I$ and $J$ are ideals of $A$
and $\theta:I\to J$ is an isomorphism. 
\end{definition}

As we saw above, 
a partial automorphism $(\theta,I,J)$ 
defines a partial morphism $\varphi:A\to\cM(J)$ 
such that $\varphi(a)b=\theta(a\theta^{-1}(b))$ 
for $a\in A$ and $b\in J$. 
If a partial morphism $\varphi:A\to\cM(B_0)$ from $A$ to $B$ 
satisfies $(\varphi(A)\cap B_0)B_0=B_0$, 
then $\varphi$ is the partial morphism determined by 
$\varphi|_I:I\to B_0$ where $I=\varphi^{-1}(B)$. 
Partial morphisms defined by partial automorphisms 
are this kind of partial morphisms. 

Following \cite{Sc1}, 
we construct $C^*$-cor\-re\-spon\-denc\-es 
from partial morphisms in the case that $A=B$. 
Let $\varphi:A\to\cM(A_0)$ be a partial morphism from $A$ to itself. 
The closed right ideal $X=A_0A$ is a Hilbert $A$-module 
with the inner product defined by $\ip{x}{y}_X=x^*y\in A$ and 
the right action defined by the multiplication. 
We have $\cK(X)\cong A_0$ and $\cL(X)\cong\cM(A_0)$ naturally. 
We define $\varphi_X:A\to \cL(X)$ 
to be the composition of $\varphi$ 
and the natural isomorphism $\cM(A_0)\to \cL(X)$. 
Thus we get a $C^*$-cor\-re\-spon\-denc\-e $X$ over $A$ 
which is denoted by $X(\varphi)$.
By definition, 
the $C^*$-cor\-re\-spon\-denc\-e $X(\varphi)$ over $A$ is non-degenerate. 
It is faithful if and only if $\varphi$ is injective, 
and full if and only if the ideal generated by $A_0$ is $A$. 
It may be reasonable to make the following definition. 

\begin{definition}
Let $\varphi:A\to\cM(A_0)$ be a partial morphism from $A$ to itself. 
We define the {\em crossed product} $A\rtimes_\varphi\N$ of $A$ by $\varphi$ 
to be the $C^*$-algebra $\cO_{X(\varphi)}$ 
of the $C^*$-cor\-re\-spon\-denc\-e $X(\varphi)$. 
\end{definition}

This definition generalizes 
crossed products by automorphisms (see \cite[Example (3)]{Pi}). 
More generally for an injective endomorphism $\varphi$ of $A$, 
its crossed product defined here 
is isomorphic to the crossed product by the endomorphism 
defined in \cite{St} 
because both of them are isomorphic to 
the augmented Cuntz-Pimsner algebra of 
the faithful $C^*$-cor\-re\-spon\-denc\-e $X(\varphi)$. 
When the partial morphism $\varphi$ of $A$ is defined 
by a partial automorphism $(\theta,I,J)$, 
then the $C^*$-cor\-re\-spon\-denc\-e $X(\varphi)$ over $A$ 
defined by $\varphi$ satisfies $J_{X(\varphi)}=I$ 
and $\varphi_{X(\varphi)}(J_{X(\varphi)})=\cK(X(\varphi))$. 
Hence the $C^*$-cor\-re\-spon\-denc\-e $X(\varphi)$ is 
a Hilbert $A$-bimodule. 
Therefore the crossed product $A\rtimes_\varphi\N$ 
of $A$ by the partial morphism $\varphi$ defined by 
a partial automorphism $(\theta,I,J)$ 
is isomorphic to 
the crossed product $A\rtimes_\theta\N$ 
of the partial automorphism $(\theta,I,J)$ 
defined in \cite[Definition 3.7]{Ex1} 
because both of them are isomorphic to 
the crossed product 
by the Hilbert $C^*$-bimodule $X(\varphi)$ 
(see Example \ref{relCP} and \cite[Example 3.2]{AEE}).

Note that a partially defined continuous map $\sigma$ 
from a locally compact space $X$ to itself 
determines a topological graph $E_\sigma=(X,U,\iota,\sigma)$ 
where $U$ is the domain of the map $\sigma$ 
and $\iota:U\to X$ is the embedding. 
The Hilbert $C^*$-bimodule $X(\varphi)$ 
of the partial morphism $\varphi$ defined from $\sigma$ 
coincides with the one obtained from 
the topological graph $E_\sigma$ as 
in Example \ref{SecTopGra}. 
Conversely every topological graph $E=(E^0,E^1,d,r)$ 
with an injective domain map $d$ 
arises in this manner. 
As explained in \cite{K1}, 
a topological graph can be considered as 
a (kind of) multi-valued continuous map 
between locally compact spaces. 
Its natural non-commutative analogue is 
a non-degenerate $C^*$-cor\-re\-spon\-denc\-e. 
Hence we can say that 
a non-degenerate $C^*$-cor\-re\-spon\-denc\-e $X$ over $A$ 
is a {\em multi-valued morphism} from $A$ to itself, 
and that the $C^*$-algebra $\cO_X$ is 
the {\em crossed product} by this multi-valued morphism.

\bibliographystyle{amsalpha}

\end{document}